\documentclass[11pt]{article}
\usepackage{cite}
\usepackage{mathrsfs}
\usepackage{amsfonts}
\usepackage{amsmath}
\usepackage{amsfonts,amssymb}
\usepackage{dsfont}
\usepackage{curves}
\usepackage{mathrsfs}
\usepackage{pifont}
\usepackage{amssymb}
\allowdisplaybreaks

\numberwithin{equation}{section}

\date{}

\begin{document}
\title{A sufficient condition for the existence of fractional $(g,f,n)$-critical covered graphs
}
\author{\small  Jie Wu\footnote{Corresponding
author. E-mail address: jskjdxwj@126.com}\\
\small  School of Economics and management\\
\small  Jiangsu University of Science and Technology\\
\small  Zhenjiang, Jiangsu 212100, China\\
}

\maketitle
\begin{abstract}
\noindent In data transmission networks, the availability of data transmission is equivalent to the existence of the fractional factor of the
corresponding graph which is generated by the network. Research on the existence of fractional factors under specific network structures can
help scientists design and construct networks with high data transmission rates. A graph $G$ is called a fractional $(g,f)$-covered graph if
for any $e\in E(G)$, $G$ admits a fractional $(g,f)$-factor covering
$e$. A graph $G$ is called a fractional $(g,f,n)$-critical covered graph if after removing any $n$ vertices of $G$, the resulting graph of
$G$ is a fractional $(g,f)$-covered graph. In this paper, we verify that if a graph $G$ of order $p$ satisfies
$p\geq\frac{(a+b-1)(a+b-2)+(a+d)n+1}{a+d}$, $\delta(G)\geq\frac{(b-d-1)p+(a+d)n+a+b+1}{a+b-1}$ and
$\delta(G)>\frac{(b-d-2)p+2\alpha(G)+(a+d)n+1}{a+b-2}$, then $G$ is a fractional $(g,f,n)$-critical covered graph, where
$g,f:V(G)\rightarrow Z^{+}$ be two functions such that $a\leq g(x)\leq f(x)-d\leq b-d$ for all $x\in V(G)$, which is a generalization of
Zhou's previous result [S. Zhou, Some new sufficient conditions for graphs to have fractional $k$-factors, International Journal of Computer
Mathematics 88(3)(2011)484--490].
\\
\begin{flushleft}
{\em Keywords:} graph; minimum degree; independence number; fractional $(g,f)$-factor; fractional $(g,f,n)$-critical covered graph.

(2010) Mathematics Subject Classification: 05C70, 68M10, 90B99
\end{flushleft}
\end{abstract}

\section{Introduction}

We research the fractional factor problem of graphs, which can be regarded as a relaxation of the well-known cardinality matching problem which
is one of the hot topics in operations research and graph theory. The fractional factor problem is widely used in the fields of network design,
scheduling and the combinatorial polyhedron. For instance, in a data transmission network we send some large data packets to several destinations
via to channels of the data transmission network, and the efficiency will be improved if the large data packets can be divided into much smaller
ones. The problem of feasible assignment of data packets can be looked as the existence of fractional flow in the network, and it can be converted
into the fractional factor problem of the network graph.

With particular reference to, a whole network can be modelled as a graph in which each vertex corresponds to a site and each edge corresponds to
a channel. The concept of a fractional $(g,f,n)$-critical covered graph reflects the feasibility of data transmission of a network. If some sites
are damaged and a channel is assigned in the process of the data transmission at the moment, the possibility of data transmission between sites
is characterized by whether the corresponding graph of the network is a fractional $(g,f,n)$-critical covered graph or not. It will imply which
structures of network can ensure the success of data transmission, and the theoretical results derived in our work can help us carry out the
network design. In the normal network, the route of data transmission is chosen by the shortest path between vertices. Several contributions on data
transmission in networks are put forward recently. Rolim et al. \cite{RRLBGSS} improved the data transmission by investigating an urban sensing problem
according to opportunistic networks. Lee, Tan and Khisti \cite{LTK} discussed the streaming data transmission on a discrete memoryless channel.
But, in the context of software definition network (SDN), the path between vertices in data transmission depends on the current network flow computation.
With minimum transmission congestion in the current moment, the transmission route is selected. In this view, the framework of data transmission problem
in SDN is equivalent to the existence of the fractional $(g,f,n)$-critical covered graph.

The graphs discussed are assumed to be simple and finite. Let $G$ be a graph. We let $V(G)$ denote the vertex set of $G$, and let $E(G)$
denote the edge set of $G$. The degree $d_G(x)$ of a vertex $x$ in $G$ is the number of edges of $G$ adjacent to $x$. The neighbour set
$N_G(x)$ of a vertex $x$ in $G$ is the set of vertices of $G$ incident with $x$. Note that $d_G(x)=|N_G(x)|$. Let $X$ be a vertex subset of $G$.
We write $G[X]$ for the subgraph of $G$ induced by $X$, and $G-X$ for the subgraph derived from $G$ by removing all the vertices in $X$ together
with the edges adjacent to vertices in $X$, namely, $G-X=G[V(G)\setminus X]$. We call $X$ independent if no two vertices in $X$ are adjacent.
We write $\delta(G)$ for the minimum degree of $G$, and $\alpha(G)$ for the independence number of $G$.

Let $g,f:V(G)\rightarrow Z^{+}$ be two functions such that $g(x)\leq f(x)$ for every $x\in V(G)$. Then a $(g,f)$-factor of $G$ is a spanning
subgraph $F$ of $G$ with $g(x)\leq d_F(x)\leq f(x)$ for all $x\in V(G)$. Let $h:E(G)\rightarrow[0,1]$ be a function satisfying
$g(x)\leq d_G^{h}(x)\leq f(x)$ for any $x\in V(G)$, where $d_G^{h}(x)=\sum\limits_{e\in E(x)}{h(e)}$ and $E(x)=\{e:e=xy\in E(G)\}$. If we write
$F_h=\{e:e\in E(G), h(e)\neq0\}$, then we say $G[F_h]$ a fractional $(g,f)$-factor of $G$ with indicator function $h$. A graph $G$ is called a
fractional $(g,f)$-covered graph if for every $e\in E(G)$, $G$ possesses a fractional $(g,f)$-factor $G[F_h]$ with $h(e)=1$. A graph $G$ is called
a fractional $(g,f,n)$-critical covered graph if after removing any $n$ vertices of $G$, the resulting graph of $G$ is a fractional $(g,f)$-covered
graph. If $g(x)=a$ and $f(x)=b$ for all $x\in V(G)$, then we call such a fractional $(g,f,n)$-critical covered graph a fractional $(a,b,n)$-critical
covered graph. A fractional $(k,k,n)$-critical covered graph is also called a fractional $(k,n)$-critical covered graph.

There exists a rich literature on the existence of factors and fractional factors in graphs. More specifically, lots of results on factors of
graphs can be discovered in \cite{AN,ET,CK,N,ZL1,ZL2,ZSL,WZi,WZo,Zd,Zs,ZBS,ZWB,Wp} and many results on fractional factors of graphs can be discovered in
\cite{LZ,LY,Zb,ZLX,WZr,WZt,K,GGC,Za1,J}. Zhou \cite{Zh} posed some independence number and minimum degree conditions for a graph admitting a fractional
$k$-factor. Yuan and Hao \cite{YH} derived a degree condition for the existence of fractional $[a,b]$-covered graphs. Zhou \cite{Za,Zha}
demonstrated some results related to fractional $(a,b,n)$-critical covered graphs.

\medskip

\noindent{\textbf{Theorem 1}} (\cite{Zh}). Let $k$ be an integer with $k\geq2$, and let $G$ be a graph of order $p$. Then $G$ contains a fractional
$k$-factor if\\
$(1)$ $k$ is even, $p>4k+1-4\sqrt{k+2}$, $\delta(G)\geq\frac{(k-1)(p+2)}{2k-1}$ and $\delta(G)>\frac{(k-2)p+2\alpha(G)-2}{2k-2}$; or\\
$(2)$ $k$ is odd, $p>4k-5$, $\delta(G)>\frac{(k-1)(p+2)}{2k-1}$ and $\delta(G)>\frac{(k-2)p+2\alpha(G)-1}{2k-2}$.

\medskip

In this paper, we are to pose a study on the existence of fractional $(g,f,n)$-critical covered graphs and verify the following theorem.
Furthermore, our result is an extension of Theorem 1.

\medskip

\noindent{\textbf{Theorem 2.}} Let $a,b,d$ and $n$ be four integers satisfying $d\geq0$, $n\geq0$, $a\geq1$ and $b-d\geq\max\{a,2\}$, let $G$
be a graph of order $p$ satisfying $p\geq\frac{(a+b-1)(a+b-2)+(a+d)n+1}{a+d}$, and let $g,f:V(G)\rightarrow Z^{+}$ be two functions such that
$a\leq g(x)\leq f(x)-d\leq b-d$ for all $x\in V(G)$. If
$$
\delta(G)\geq\frac{(b-d-1)p+(a+d)n+a+b+1}{a+b-1}
$$
and
$$
\delta(G)>\frac{(b-d-2)p+2\alpha(G)+(a+d)n+1}{a+b-2},
$$
then $G$ is a fractional $(g,f,n)$-critical covered graph.

\medskip

Corollary 1 is easily derived by setting $d=0$ in Theorem 2.

\medskip

\noindent{\textbf{Corollary 1.}} Let $a,b$ and $n$ be three integers satisfying $n\geq0$, $a\geq1$ and $b\geq\max\{a,2\}$, let $G$
be a graph of order $p$ satisfying $p\geq\frac{(a+b-1)(a+b-2)+an+1}{a}$, and let $g,f:V(G)\rightarrow Z^{+}$ be two functions such that
$a\leq g(x)\leq f(x)\leq b$ for all $x\in V(G)$. If
$$
\delta(G)\geq\frac{(b-1)p+an+a+b+1}{a+b-1}
$$
and
$$
\delta(G)>\frac{(b-2)p+2\alpha(G)+an+1}{a+b-2},
$$
then $G$ is a fractional $(g,f,n)$-critical covered graph.

\medskip

Corollary 2 is easily derived by setting $g(x)\equiv f(x)$ for all $x\in V(G)$ in Corollary 1.

\medskip

\noindent{\textbf{Corollary 2.}} Let $a,b$ and $n$ be three integers satisfying $n\geq0$, $a\geq1$ and $b\geq\max\{a,2\}$, let $G$
be a graph of order $p$ satisfying $p\geq\frac{(a+b-1)(a+b-2)+an+1}{a}$, and let $f:V(G)\rightarrow Z^{+}$ be a function such that
$a\leq f(x)\leq b$ for all $x\in V(G)$. If
$$
\delta(G)\geq\frac{(b-1)p+an+a+b+1}{a+b-1}
$$
and
$$
\delta(G)>\frac{(b-2)p+2\alpha(G)+an+1}{a+b-2},
$$
then $G$ is a fractional $(f,n)$-critical covered graph.

\medskip

Corollary 3 is easily derived by setting $n=0$ in Theorem 2.

\medskip

\noindent{\textbf{Corollary 3.}} Let $a,b$ and $d$ be three integers satisfying $d\geq0$, $a\geq1$ and $b-d\geq\max\{a,2\}$, let $G$
be a graph of order $p$ satisfying $p\geq\frac{(a+b-1)(a+b-2)+1}{a+d}$, and let $g,f:V(G)\rightarrow Z^{+}$ be two functions such that
$a\leq g(x)\leq f(x)-d\leq b-d$ for all $x\in V(G)$. If
$$
\delta(G)\geq\frac{(b-d-1)p+a+b+1}{a+b-1}
$$
and
$$
\delta(G)>\frac{(b-d-2)p+2\alpha(G)+1}{a+b-2},
$$
then $G$ is a fractional $(g,f)$-covered graph.

\section{The proof of Theorem 2}

For the proof of Theorem 2, we shall utilize the following theorem which was verified by Li, Yan and Zhang \cite{LYZ}.

\medskip

\noindent{\textbf{Theorem 3}} (\cite{LYZ}). Let $G$ be a graph, and let $g,f:V(G)\rightarrow Z^{+}$ be two functions with $g(x)\leq f(x)$ for every
$x\in V(G)$. Then $G$ is a fractional $(g,f)$-covered graph if and only if
$$
\gamma_G(X,Y)=f(X)+d_{G-X}(Y)-g(Y)\geq\varepsilon(X)
$$
for any $X\subseteq V(G)$, where $Y=\{x:x\in V(G)\setminus X, d_{G-X}(x)\leq g(x)\}$ and $\varepsilon(X)$ is defined by
\[
 \varepsilon(X)=\left\{
\begin{array}{ll}
2,&if \ X \ is \ not \ independent,\\
1,&if \ X \ is \ independent \ and \ there \ is \ an \ edge \ joining\\
&X \ and \ V(G)\setminus(X\cup Y), \ or \ there \ is \ an \ edge \ e=xy\\
&joining \ X \ and \ Y \ such \ that \ d_{G-X}(y)=g(y) \ for\\
&y\in Y,\\
0,&otherwise.\\
\end{array}
\right.
\]

\medskip

\noindent{\it Proof of Theorem 2.} \  Let $N\subseteq V(G)$ with $|N|=n$, and put $H=G-N$. To verify Theorem 2, it suffices to prove that $H$
is a fractional $(g,f)$-covered graph. Suppose, to the contrary, that $H$ is not a fractional $(g,f)$-covered graph. Then by Theorem 3, we admit
$$
\gamma_H(X,Y)=f(X)+d_{H-X}(Y)-g(Y)\leq\varepsilon(X)-1\eqno(1)
$$
for some subset $X\subseteq V(H)$, where $Y=\{x:x\in V(H)\setminus X, d_{H-X}(x)\leq g(x)\}$.

If $Y=\emptyset$, then from (1) and $\varepsilon(X)\leq|X|$, $\varepsilon(X)-1\geq\gamma_H(X,\emptyset)=f(X)\geq(a+d)|X|\geq|X|\geq\varepsilon(X)$,
a contradiction. Therefore, $Y\neq\emptyset$. Next, we may define
$$
h=\min\{d_{H-X}(x):x\in Y\}
$$
and choose $x_1\in Y$ satisfying $d_{H-X}(x_1)=h$. Obviously, $0\leq h\leq b-d$. Furthermore, $\delta(H)\leq d_H(x_1)\leq d_{H-X}(x_1)+|X|=h+|X|$,
namely,
$$
|X|\geq\delta(H)-h.\eqno(2)
$$

In what follows, We shall consider two cases by the value of $h$.

\noindent{\bf Case 1.} \ $1\leq h\leq b-d$.

In light of (1), (2), $H=G-N$ with $|N|=n$, $\varepsilon(X)\leq2$ and $|X|+|Y|+n\leq p$, we deduce
\begin{eqnarray*}
1&\geq&\varepsilon(X)-1\geq\gamma_H(X,Y)=f(X)+d_{H-X}(Y)-g(Y)\\
&\geq&(a+d)|X|+h|Y|-(b-d)|Y|\\
&=&(a+d)|X|-(b-d-h)|Y|\\
&\geq&(a+d)|X|-(b-d-h)(p-n-|X|)\\
&=&(a+b-h)|X|-(b-d-h)p+(b-d-h)n\\
&\geq&(a+b-h)(\delta(H)-h)-(b-d-h)p+(b-d-h)n\\
&=&(a+b-h)(\delta(G-N)-h)-(b-d-h)p+(b-d-h)n\\
&\geq&(a+b-h)(\delta(G)-n-h)-(b-d-h)p+(b-d-h)n\\
&=&(a+b-h)\delta(G)-(a+b-h)h-(b-d-h)p-(a+d)n\\
&\geq&(a+b-h)\cdot\frac{(b-d-1)p+(a+d)n+a+b+1}{a+b-1}-(a+b-h)h\\
&&-(b-d-h)p-(a+d)n.
\end{eqnarray*}
Multiplying the above inequality by $(a+b-1)$ and rearranging,
$$
0\geq(h-1)((a+d)p-(a+b-1)(a+b-h)-(a+d)n-2)+a+b-1.\eqno(3)
$$

If $h=1$, then (3) is impossible by $a\geq1$ and $b\geq\max\{a,2\}+d$. Next, we consider $2\leq h\leq b-d$.

Note that $p\geq\frac{(a+b-1)(a+b-2)+(a+d)n+1}{a+d}$ and $2\leq h\leq b-d$. Thus, we admit
\begin{eqnarray*}
&&(a+d)p-(a+b-1)(a+b-h)-(a+d)n-2\\
&\geq&(a+b-1)(a+b-2)+(a+d)n+1-(a+b-1)(a+b-h)-(a+d)n-2\\
&=&(a+b-1)(h-2)-1\geq-1,
\end{eqnarray*}
and so,
$$
(h-1)((a+d)p-(a+b-1)(a+b-h)-(a+d)n-2)\geq-(h-1).\eqno(4)
$$

It follows from (3), (4) and $2\leq h\leq b-d$ that
\begin{eqnarray*}
0&\geq&(h-1)((a+d)p-(a+b-1)(a+b-h)-(a+d)n-2)+a+b-1\\
&\geq&-(h-1)+a+b-1=a+b-h\geq a+b-(b-d)=a+d\\
&\geq&a\geq1,
\end{eqnarray*}
which is a contradiction.

\noindent{\bf Case 2.} \ $h=0$.

Let $Q=\{x\in Y: d_{H-X}(x)=0\}$, $W=\{x\in Y: d_{H-X}(x)=1\}$, $W_1=\{x\in W: N_{H-X}(x)\subseteq Y\}$ and $W_2=W-W_1$. Then the graph induced
by $W_1$ in $H-X$ admits maximum degree at most 1. Let $A$ be a maximum independent set of this graph. Obviously, $|A|\geq\frac{1}{2}|W_1|$.
According to our definitions, $Q\cup A\cup W_2$ is an independent set of $H$. Therefore, we infer
$$
\alpha(H)\geq|Q|+|A|+|W_2|\geq|Q|+\frac{1}{2}|W_1|+\frac{1}{2}|W_2|=|Q|+\frac{1}{2}|W|.\eqno(5)
$$

We easily see that
$$
\alpha(G)\geq\alpha(H)\eqno(6)
$$
by $H=G-N$ with $N\subseteq V(G)$ and $|N|=n$. Using (5) and (6), we deduce
$$
\alpha(G)\geq\alpha(H)\geq|Q|+\frac{1}{2}|W|.\eqno(7)
$$

In light of (1), (2), (7), $H=G-N$ with $|N|=n$, $|X|+|Y|+n\leq p$, $\varepsilon(X)\leq2$
and $h=0$, we infer
\begin{eqnarray*}
1&\geq&\varepsilon(X)-1\geq\gamma_H(X,Y)=f(X)+d_{H-X}(Y)-g(Y)\\
&\geq&(a+d)|X|+d_{H-X}(Y\setminus(Q\cup W))+|W|-(b-d)|Y|\\
&\geq&(a+d)|X|+2|Y\setminus(Q\cup W)|+|W|-(b-d)|Y|\\
&=&(a+d)|X|-2\Big(|Q|+\frac{1}{2}|W|\Big)-(b-d-2)|Y|\\
&\geq&(a+d)|X|-2\alpha(G)-(b-d-2)(p-n-|X|)\\
&=&(a+b-2)|X|-2\alpha(G)-(b-d-2)(p-n)\\
&\geq&(a+b-2)(\delta(H)-h)-2\alpha(G)-(b-d-2)(p-n)\\
&=&(a+b-2)\delta(H)-2\alpha(G)-(b-d-2)(p-n)\\
&\geq&(a+b-2)(\delta(G)-n)-2\alpha(G)-(b-d-2)(p-n)\\
&=&(a+b-2)\delta(G)-2\alpha(G)-(b-d-2)p-(a+d)n,
\end{eqnarray*}
which implies
$$
\delta(G)\leq\frac{(b-d-2)p+2\alpha(G)+(a+d)n+1}{a+b-2},
$$
which contradicts that $\delta(G)>\frac{(b-d-2)p+2\alpha(G)+(a+d)n+1}{a+b-2}$. This completes the proof of Theorem 2. \hfill $\Box$

\section{Remarks}

\noindent{\bf Remark 1.} \ In what follows, we explain that the condition $\delta(G)\geq\frac{(b-d-1)p+(a+d)n+a+b+1}{a+b-1}$ in Theorem 2
is best possible, that is to say, we cannot replaced $\delta(G)\geq\frac{(b-d-1)p+(a+d)n+a+b+1}{a+b-1}$ by
$\delta(G)\geq\frac{(b-d-1)p+(a+d)n+a+b}{a+b-1}$ in Theorem 2.

Let $b=a+d$, $g(x)\equiv b-d$ and $f(x)\equiv a+d$, and let $t$ be a sufficiently large positive integer such that $\frac{2(b-d-1)t+1}{a+d}$
is an integer. We let $G=K_{\frac{2(b-d-1)t+1}{a+d}+n}\bigvee(tK_2)$ be a graph, and set $p=|V(G)|=\frac{2(b-d-1)t+1}{a+d}+n+2t=\frac{2t(a+b-1)+1}{a+d}+n$,
where $\bigvee$ means ``join". We easily deduce that $\alpha(G)=t$, $\delta(G)=\frac{2(b-d-1)t+1}{a+d}+n+1=\frac{(b-d-1)p+(a+d)n+a+b}{a+b-1}$
and $\delta(G)=\frac{2(b-d-1)t+1}{a+d}+n+1>\frac{(b-d-2)p+2\alpha(G)+(a+d)n+1}{a+b-2}$.

For any $N\subseteq V(K_{\frac{2(b-d-1)t+1}{a+d}+n})$ with $|N|=n$, we let $H=G-N$. Write $X=V(K_{\frac{2(b-d-1)t+1}{a+d}+n})\setminus N$
and $Y=V(tK_2)$. Clearly, $\varepsilon(X)=2$. Thus, we infer
\begin{eqnarray*}
\gamma_H(X,Y)&=&f(X)+d_{H-X}(Y)-g(Y)\\
&=&(a+d)|X|+d_{H-X}(Y)-(b-d)|Y|\\
&=&(a+d)\cdot\frac{2(b-d-1)t+1}{a+d}+2t-(b-d)\cdot2t\\
&=&1<2=\varepsilon(X).
\end{eqnarray*}
Hence, $H$ is not a fractional $(g,f)$-covered graph by Theorem 3, which implies that $G$ is not a fractional $(g,f,n)$-critical covered graph.

\medskip

\noindent{\bf Remark 2.} \ Next, we claim that the condition $\delta(G)>\frac{(b-d-2)p+2\alpha(G)+(a+d)n+1}{a+b-2}$ in Theorem 2
is sharp, that is to say, we cannot replaced $\delta(G)>\frac{(b-d-2)p+2\alpha(G)+(a+d)n+1}{a+b-2}$ by
$\delta(G)\geq\frac{(b-d-2)p+2\alpha(G)+(a+d)n+1}{a+b-2}$ in Theorem 2.

Let $b=a+d$, $g(x)\equiv b-d$ and $f(x)\equiv a+d$, and let $t$ be a sufficiently large positive integer such that $\frac{2(b-d-1)t+1}{a+d}$
is an integer. We let $G=K_{\frac{2(b-d-1)t+(b-d)(a+b)+1}{a+d}+n}\bigvee((a+b)K_1\cup(tK_2))$ be a graph, and set $p=|V(G)|=\frac{2(b-d-1)t+(b-d)(a+b)+1}{a+d}+n+a+b+2t=\frac{2t(a+b-1)+(a+b)^{2}+1}{a+d}+n$. We easily infer that $\alpha(G)=a+b+t$, $\delta(G)=\frac{2(b-d-1)t+(b-d)(a+b)+1}{a+d}+n=\frac{(b-d-2)p+2\alpha(G)+(a+d)n+1}{a+b-2}$
and $\delta(G)=\frac{2(b-d-1)t+(b-d)(a+b)+1}{a+d}+n=\frac{(b-d-1)p+(a+d)n+a+b+1}{a+b-1}$.

For any $N\subseteq K_{\frac{2(b-d-1)t+(b-d)(a+b)+1}{a+d}+n})$ with $|N|=n$, we let $H=G-N$. Write $X=V(K_{\frac{2(b-d-1)t+(b-d)(a+b)+1}{a+d}+n})\setminus N$
and $Y=V((a+b)K_1\cup(tK_2))$. Obviously, $\varepsilon(X)=2$. Thus, we derive
\begin{eqnarray*}
\gamma_H(X,Y)&=&f(X)+d_{H-X}(Y)-g(Y)\\
&=&(a+d)|X|+d_{H-X}(Y)-(b-d)|Y|\\
&=&(a+d)\cdot\frac{2(b-d-1)t+(b-d)(a+b)+1}{a+d}+2t-(b-d)\cdot(a+b+2t)\\
&=&1<2=\varepsilon(X).
\end{eqnarray*}
Therefore, $H$ is not a fractional $(g,f)$-covered graph by Theorem 3, namely, $G$ is not a fractional $(g,f,n)$-critical covered graph.


\end{document}